\documentclass[a4paper,11pt]{article}
\usepackage{amsmath, amsfonts, amscd, amssymb, amsthm, enumerate, xypic}

\def\car{\text{char}}
\def\Mat{\text{M}}
\def\GL{\text{GL}}
\def\SL{\text{SL}}

\def\card{\#\,}
\def\defterm{\textbf}

\newcommand{\tr}{\operatorname{tr}}
\newcommand{\Ker}{\operatorname{Ker}}

\renewcommand{\setminus}{\smallsetminus}

\def\F{\mathbb{F}}
\def\K{\mathbb{K}}

\def\N{\mathbb{N}}

\def\calA{\mathcal{A}}

\def\calE{\mathcal{E}}

\def\lcro{\mathopen{[\![}}
\def\rcro{\mathclose{]\!]}}

\theoremstyle{definition}
\newtheorem{Def}{Definition}
\newtheorem{Not}[Def]{Notation}

\theoremstyle{plain}
\newtheorem{theo}{Theorem}

\newtheorem{cor}[theo]{Corollary}
\newtheorem{lemme}[theo]{Lemma}

\theoremstyle{remark}
\newtheorem{Rems}{Remarks}
\newtheorem{Rem}[Rems]{Remark}

\title{On decomposing any matrix as a linear combination of three idempotents}
\author{Cl\'ement de Seguins Pazzis
\footnote{Professor of Mathematics at Lyc\'ee Priv\'e Sainte-Genevi\`eve, 2, rue
de l'\'Ecole des Postes, 78029 Versailles Cedex, FRANCE.}
\footnote{e-mail address: dsp.prof@gmail.com}}
\date{\today}

\begin{document}

\maketitle

\begin{abstract}
In a recent article, we gave a full characterization of matrices that can be
decomposed as linear combinations of two idempotents with prescribed coefficients.
In this one, we use those results to improve on a recent theorem of V. Rabanovich: we establish that every square matrix is
a linear combination of three idempotents (for an arbitrary coefficient field rather than just one
of characteristic $0$).
\end{abstract}



\section{Introduction}

In this article, $\K$ will denote an arbitrary field, $\car (\K)$ its characteristic,
and $n$ a positive integer. We choose an algebraic closure $\overline{\K}$ of $\K$.
We will use the French convention for the set of integers:
$\N$ will denote the set of non-negative integers, and $\N^*$ the one of positive integers.

\vskip 2mm
\noindent
An idempotent matrix of $\Mat_n(\K)$ is a matrix $P$
verifying $P^2=P$, i.e. idempotent matrices represent projectors in finite dimensional
vector spaces. Of course, any matrix similar to an idempotent is itself an idempotent.

\vskip 2mm
\noindent Our main topic of interest is determining
the smallest integer $\ell_n(\K)$ such that any matrix of $\Mat_n(\K)$ can be decomposed
into a linear combination (LC) of $\ell_n(\K)$ idempotents.

\noindent Our main results are summed up in the following theorem:

\begin{theo}[Main theorem]\label{maindSP}
Any matrix of $M_n(\K)$ is a linear combination of $3$ idempotents. \\
More precisely, equality $\ell_n(\K)=3$ holds save for the following special cases:
\begin{enumerate}[(a)]
\item If $n=1$, then $\ell_n(\K)=1$;
\item If $n=2$ and $\card \K>2$, then $\ell_n(\K)=2$;
\item If $n=3$ and every polynomial of degree $3$ in $\K[X]$
has a root in $\K$, then $\ell_n(\K)=2$.
\end{enumerate}
\end{theo}

Inequality $\ell_n(\K) \leq 3$ was already known prior to this paper
for a field of characteristic $0$ (see \cite{Rabanovich}) with a more elementary proof
that cannot be generalized to an arbitrary field.

\begin{Rem}[A trivial but nevertheless useful remark]
Since the zero matrix is an idempotent, any matrix that is a linear combination of $p$ idempotents is also a linear
combination of $k$ idempotents for every integer $k \geq p$.
\end{Rem}

\noindent The rest of the paper is laid out as follows:
\begin{enumerate}[(1)]
\item We will start by reviewing some characterizations of linear combinations of two idempotents
 that were featured in \cite{dSPidem2}.
\item These results will then be used to give a lower bound for $\ell_n(\K)$.
\item Proving that $\ell_n(\K)\leq 3$ is much more demanding and will require
subtle manipulations of cyclic matrices and rational canonical forms
(see \cite{Lev} for similar constructions in a different context). Therefore, section 5 features
a review of cyclic matrices. Finally, section 6 consists of the proof
that every square matrix is a linear combination of three idempotents.
Given $M \in \Mat_n(\K)$, our basic strategy will be to find an idempotent $P$ and a scalar $a$ such that
$M-a.P$ is a linear combination of two idempotents.
\end{enumerate}

\section{Additional notations}

Given a list $(A_1,\dots,A_p)$ of square matrices, we will let
$$D(A_1,\dots,A_p):=\begin{bmatrix}
A_1 & 0 & & 0 \\
0 & A_2 & & \vdots \\
\vdots & & \ddots & \\
0 & \dots & &  A_p
\end{bmatrix}
$$
denote the block-diagonal matrix with diagonal blocks $A_1$, \dots, $A_p$.

\vskip 2mm
\noindent Similarity of two matrices $A$ and $B$ of $\Mat_n(\K)$ will be written $A\sim B$.

\vskip 2mm
\noindent The characteristic polynomial of a matrix $M$ will be denoted by $\chi_M$, its trace by $\tr M$.

\vskip 2mm
\noindent
Let $P=X^n-\underset{k=0}{\overset{n-1}{\sum}}a_kX^k \in \K[X]$ be a monic polynomial with degree $n$.
Its \defterm{companion matrix} is
$$C(P):=\begin{bmatrix}
0 &   & & 0 & a_0 \\
1 & 0 & &   & a_1 \\
0 & \ddots & \ddots & & \vdots \\
\vdots & & & 0 & a_{n-2} \\
0 & & &  1 & a_{n-1}
\end{bmatrix}.$$
Its characteristic polynomial is precisely $P$, and so is its minimal polynomial.
We will set $\tr P:=\tr C(P)=a_{n-1}$.

\vskip 2mm
\noindent
Let $H_{n,p}$ denote the elementary matrix
$\begin{bmatrix}
0 & \cdots & 0 &  1 \\
\vdots & & 0 & 0 \\
0 & \cdots & 0 & 0
\end{bmatrix} \in \Mat_{n,p}(\K)$ with only one non-zero coefficient located on the first row and
$p$-th column. \\
For $k \in \N^*$, we set
$$F_k:=D(0,\dots,0,1)\in \Mat_k(\K).$$

\section{On linear combinations of two idempotents with prescribed coefficients}

In order to prove our theorem, we will make extensive use
of the results featured in \cite{dSPidem2}, so reviewing them is necessary.

\begin{Def}
Let $\calA$ be a $\K$-algebra and
$(\alpha_1,\dots,\alpha_n)\in (\K^*)^n$.
An element $x \in \calA$ will be called an
\textbf{$(\alpha_1,\dots,\alpha_n)$-composite} when there are idempotents $p_1,\dots,p_n$
such that $x=\underset{k=1}{\overset{n}{\sum}}\alpha_k.p_k$. \\
\end{Def}

\begin{Not}
When $A$ is a matrix of $\Mat_n(\K)$, $\lambda \in \overline{\K}$ and $k \in \N^*$, we set
$$n_k(A,\lambda):=\dim \Ker (A-\lambda.I_n)^k-\dim \Ker (A-\lambda.I_n)^{k-1},$$
i.e. $n_k(A,\lambda)$ is the number of blocks of size greater or equal to $k$
for the eigenvalue $\lambda$ in the Jordan reduction of $A$
(in particular, it is zero when $\lambda$ is not an eigenvalue of $A$).
We also denote by $j_k(A,\lambda)$ the number of blocks of size $k$ for the eigenvalue $\lambda$ in the Jordan reduction of $A$.
\end{Not}

\begin{Def}
Two sequences $(u_k)_{k \geq 1}$ and $(v_k)_{k \geq 1}$ are said to be \defterm{intertwined}
when:
$$\forall k \in \N^*, \; v_k \leq u_{k+1} \quad \text{and} \quad u_k \leq v_{k+1.}$$
\end{Def}

With that in mind, the problem of determining whether a particular matrix $A \in \Mat_n(\K)$
is an $(\alpha,\beta)$-composite is completely answered by the following theorems:

\begin{theo}\label{HPdiff}
Assume $\car(\K) \neq 2$ and let $A \in \Mat_n(\K)$.
Then $A$ is an $(\alpha,-\alpha)$-composite iff all the following conditions hold:
\begin{enumerate}[(i)]
\item The sequences $(n_k(A,\alpha))_{k \geq 1}$ and $(n_k(A,-\alpha))_{k \geq 1}$ are intertwined.
\item $\forall \lambda \in \overline{\K} \setminus \{0,\alpha,-\alpha\}, \; \forall k \in \N^*, \; j_k(A,\lambda)=j_k(A,-\lambda)$.
\end{enumerate}
\end{theo}

\begin{theo}\label{HPsum}
Assume $\car(\K) \neq 2$, and let $A \in \Mat_n(\K)$.
Then $A$ is an $(\alpha,\alpha)$-composite iff all the following conditions hold:
\begin{enumerate}[(i)]
\item The sequences $(n_k(A,0))_{k \geq 1}$ and $(n_k(A,2\,\alpha))_{k \geq 1}$
are intertwined.
\item $\forall \lambda \in \overline{\K} \setminus \{0,\alpha,2\alpha\}, \; \forall k \in \N^*, \;
j_k(A,\lambda)=j_k(A,2\alpha-\lambda)$.
\end{enumerate}
\end{theo}

\begin{theo}\label{HPsumcar2}
Assume $\car(\K)=2$ and let $A \in \Mat_n(\K)$.
Then $A$ is an $(\alpha,-\alpha)$-composite iff
for every $\lambda \in \overline{\K} \setminus \{0,\alpha\}$,
all blocks in the Jordan reduction of $A$ with respect to $\lambda$ have an even
size.
\end{theo}

\begin{theo}\label{dSP2LC}
Let $A \in \Mat_n(\K)$ and $(\alpha,\beta)\in (\K^*)^2$ such that
$\alpha \neq \pm \beta$. Then $A$ is an $(\alpha,\beta)$-composite iff all the following conditions hold:
\begin{enumerate}[(i)]
\item The sequences $(n_k(A,0))_{k \geq 1}$ and $(n_k(A,\alpha+\beta))_{k \geq 1}$
are intertwined.
\item The sequences $(n_k(A,\alpha))_{k \geq 1}$ and $(n_k(A,\beta))_{k \geq 1}$ are intertwined.
\item $\forall \lambda \in \overline{\K} \setminus \{0,\alpha,\beta,\alpha+\beta\}, \; \forall
k \in \N^*, \; j_k(A,\lambda)=j_k(A,\alpha+\beta-\lambda)$.
\item If in addition $\car(\K)\neq 2$, then
$\forall k \in \N^*, \; j_{2k+1}\bigl(A,\frac{\alpha+\beta}{2}\bigr)=0$.
\end{enumerate}
\end{theo}

\vskip 4mm
\noindent These theorems have the following easy consequences, which we will use in the next sections:

\begin{cor}\label{dim2}
Let $A \in \Mat_2(\K)$ be non-scalar with trace $t$, and let $(a,b)\in (\K^*)^2$ such that
$a+b=t$. Then $A$ is an $(a,b)$-composite.
\end{cor}

\begin{proof}[Proof using the previous theorems] ${}$
\begin{itemize}
\item If $A$ has two different eigenvalues $c$ and $d$ in $\overline{\K}$, then
$c=a+b-d$ and these eigenvalues have multiplicity $1$ therefore, using all the previous theorems, we see that
$A$ is an $(a,b)$-composite.
\item Assume now $A$ has only one eigenvalue $\lambda$. Then $a+b=2\lambda$ and
the Jordan block corresponding to $\lambda$ is even-sized, so theorems \ref{HPsumcar2} and \ref{dSP2LC}
show that $A$ is an $(a,b)$-composite.
\end{itemize}
\end{proof}

\noindent See also \cite{Rabanovich} for a very elementary proof.

\begin{cor}\label{nilpotent}
Every nilpotent matrix is a $(1,-1)$-composite, and more generally an $(\alpha,-\alpha)$-composite for every
$\alpha \in \K^*$. \\
If $\car(\K)=2$, then every unipotent\footnote{A unipotent matrix is one of
the form $I_n+N$ where $N$ is nilpotent.}
matrix is a $(1,1)$-composite.
\end{cor}

\begin{cor}\label{alphaneqbeta}
Let $\alpha \in \K^*$ and $\beta \in \K^*$ such that $\alpha \neq \beta$. Then, for every
$n \in \N^*$, the companion matrices $C\bigl((X-\alpha)^n(X-\beta)^n\bigr)$,
$C\bigl((X-\alpha)^{n+1}(X-\beta)^n\bigr)$ and $C\bigl((X-\alpha)^n(X-\beta)^{n+1}\bigr)$
are all $(\alpha,\beta)$-composites.
\end{cor}

\begin{cor}[When a diagonal matrix is an $(\alpha,\beta)$-composite]\label{diagcor} ${}$ \\
Let $A=D(a_1,\dots,a_n)$ be a diagonal matrix, and $(\alpha,\beta)\in (\K^*)^2$. \\
For $\lambda \in \K$, set $n_\lambda:=\#\{k \in \lcro 1,n\rcro : \; a_k=\lambda\}$.
\begin{enumerate}[(i)]
\item If $\car(\K)=2$, then $A$ is an $(\alpha,\alpha)$-composite iff
$a_k \in \{0,\alpha\}$ for all $k \in \lcro 1,n\rcro$.
\item If $\car(\K)\neq 2$, then $A$ is an $(\alpha,-\alpha)$-composite
iff $n_\lambda=n_{-\lambda}$ for all $\lambda \in \K \setminus \{0,\alpha,-\alpha\}$,
\item If $\car(\K) \neq 2$, then $A$ is an $(\alpha,\alpha)$-composite
iff $n_\lambda=n_{2\alpha-\lambda}$ for all $\lambda \in \K \setminus \{0,\alpha,2\alpha\}$.
\item If $\car(\K)=2$ and $\alpha \neq \beta$, then $A$ is an $(\alpha,\beta)$-composite
iff $n_\lambda=n_{\alpha+\beta-\lambda}$ for all $\lambda \in \K \setminus \{0,\alpha,\beta,\alpha+\beta\}$.
\item If $\car(\K)\neq 2$ and $\alpha \neq \beta$, then $A$ is an $(\alpha,\beta)$-composite
iff $n_{(\alpha+\beta)/2}=0$ and $n_{\lambda}=n_{\alpha+\beta-\lambda}$
for every $\lambda \in \K \setminus \bigl\{0,\alpha,\beta,\alpha+\beta,\frac{\alpha+\beta}{2}\bigr\}$.
\end{enumerate}
\end{cor}

Finally, the following corollary will be useful in some cases:

\begin{cor}\label{evenout}
Let $A \in \Mat_n(\K)$ and assume $A$ is an $(\alpha,\beta)$-composite for some $(\alpha,\beta)\in (\K^*)^2$.
Then the total multiplicity of the eigenvalues of $A$ which do not belong to $\{0,\alpha,\beta,\alpha+\beta\}$
is an even number. The total multiplicity of the eigenvalues which do not belong to $\K$ is also even.
\end{cor}

\section{A lower bound for $\ell_n(\K)$}

Here we want to prove the ``lower bound" part of our main theorem.
The case $n=1$ is trivial, so we immediately move on to the case $n \geq 2$.
A non-zero nilpotent matrix of $\Mat_n(\K)$ is not the product of an idempotent by a
scalar, thus $\ell_n(\K) \geq 2$.
\begin{enumerate}
\item Assume $n=2$ and $\card \K >3$.
If $A$ is scalar (i.e. a multiple of $I_2$), then it is a $(1,0)$-composite.
Assume $A$ is not scalar. Since $\card \K>3$, the set
$\{\tr A-a \mid a\in \K^*\}$ has at least two elements, hence a non-zero element $\alpha$,
so Corollary \ref{dim2} shows that $A$ is an $(\alpha,\tr A-\alpha)$-composite.
This proves $\ell_n(\K)=2$.

\vskip 2mm
\item Assume $n=2$ and $\K=\F_2$. \\
Then the matrix $A=\begin{bmatrix}
0 & 1 \\
1 & 1
\end{bmatrix}$ is not a linear combination of two idempotents. Indeed, if it were,
it would be a sum of two idempotents (since it is not an idempotent itself),
but this is not the case since $A$ has two distinct eigenvalues in $\overline{\F_2}\setminus \F_2$
with multiplicity $1$, hence with odd-sized Jordan blocks, in contradiction with Theorem \ref{HPsumcar2}.
This proves $\ell_2(\F_2) \geq 3$.

\vskip 2mm
\item Assume $n=3$ and there is an irreducible polynomial $P \in \K[X]$ of degree $3$. Without loss of generality, we can assume $P=X^3-aX^2-bX-c$ for some
$(a,b,c)\in \K^3$. We claim that the companion matrix
$A=\begin{bmatrix}
0 & 0 & c \\
1 & 0 & b \\
0 & 1 & a
\end{bmatrix}$ is not a linear combination of two idempotents:
since $A$ has no eigenvalue in $\K$, it is not the product of an idempotent by a scalar;
it is neither an $(\alpha,\beta)$-composite for some $(\alpha,\beta)\in (\K^*)^2$ because it is
odd-sized and has no eigenvalue in $\K$ (see Corollary \ref{evenout}).
This shows $\ell_n(\K) \geq 3$.

\vskip 2mm
\item Assume $n=3$ and every polynomial $P \in \K[X]$ of degree $3$ has a root in~$\K$. \\
As a consequence, the field $\K$ is infinite (recall that when $\K$ is finite, there
exists, for every $k \in \N^*$, an irreducible polynomial of degree $k$ in $\K[X]$).
We then claim that every matrix of $\Mat_3(\K)$ is a linear combination of two idempotents. \\
Let $A \in \Mat_3(\K)$. Leaving the trivial cases aside, we can assume
$A$ has more than one eigenvalue in $\overline{\K}$, so reduction to a canonical form shows, combined
with the assumption on roots of polynomials of degree $3$, that
$A$ is similar to $\begin{bmatrix}
\lambda & 0 & 0 \\
0 & 0 & b \\
0 & 1 & a
\end{bmatrix}$ for some triple $(\lambda,a,b)\in \K^3$. \\
If $\lambda =0$, then the previous cases show that $A$ is an LC of two idempotents. \\
Assume now $\lambda \neq 0$.
If $\lambda \neq a$, then Corollary \ref{dim2} shows that the block matrix
$\begin{bmatrix}
0 & b \\
1 & a
\end{bmatrix}$ is a $(\lambda,a-\lambda)$-composite, hence $A$ is also a $(\lambda,a-\lambda)$-composite. \\
If $\lambda=a$, then we can find a pair $(b,c)\in (\K^*)^2$ such that $a=b+c$,
and again, since $\lambda=b+c$, Corollary \ref{dim2} shows that $A$ is a $(b,c)$-composite.
In any case, we have proven that $A$ is an LC of two idempotents. We conclude that
$\ell_3(\K)=2$.

\vskip 2mm
\item Assume finally $n \geq 4$. We wish to prove then that $\ell_n(\K) \geq 3$. \\
If $\K$ is finite, then we can find a monic polynomial $P=X^3-aX^2-bX-c$ of degree $3$
with no root in $\K$, and the same line of reasoning as in point 3 shows that the matrix
$A=\begin{bmatrix}
0_{n-3} & 0 & 0 & 0 \\
0 & 0 & 0 & c \\
0 & 1 & 0 & b \\
0 & 0 & 1 & a
\end{bmatrix}$ is not an LC of two idempotents.
Assume now $\K$ is infinite, and choose arbitrary elements
$a_1,a_2,a_3,a_4$ in $\K$. Assume furthermore that:
\begin{enumerate}[(i)]
\item $a_i \neq \pm a_j$ for all distinct $i$ and $j$;
\item $a_i \neq a_j+a_k$ for all $i,j$ and $k$ (distinct or not);
\item $a_i+a_j \neq a_k+a_l$ for all distinct $i,j,k,l$.
\end{enumerate}
Condition (ii) in the case $i=j=k$ shows that the $a_i$'s are non-zero, and
condition (i) shows that the $a_i$'s are pairwise distinct.
We wish to prove that the diagonal matrix
$A=D(a_1,a_2,a_3,a_4,0,\dots,0)$ of $\Mat_n(\K)$ is not an LC of two idempotents. \\
In doing so, we will use Corollary \ref{diagcor} repeatedly.
By a \emph{reductio ad absurdum}, let us assume
$A$ is an $(\alpha,\beta)$-composite for some $(\alpha,\beta)\in (\K^*)^2$ (since clearly it is not a scalar multiple of an idempotent).
\begin{itemize}
\item If $\alpha=-\beta$ and $\car (\K) \neq 2$, then some $a_i$ is different from $\alpha$, $0$ and $-\alpha$, so
case (ii) in Corollary \ref{diagcor} shows that $-a_i$ should be another eigenvalue of $A$, which is forbidden by condition (i).
\item Assume $\alpha=\beta$ and $\car(\K) \neq 2$:
then condition (ii) ensures that at most one of the $a_i$'s
belongs to $\{\alpha,2\,\alpha\}$, so, using again Corollary \ref{diagcor}, we see that
\emph{none} of the $a_i$'s belongs to $\{\alpha,2\,\alpha\}$; case (iii) in Corollary \ref{diagcor} then shows
that there is a permutation $\sigma$ of $\{1,2,3,4\}$ such that $a_{\sigma(2)}=2\,\alpha-a_{\sigma(1)}$
and $a_{\sigma(4)}=2\,\alpha-a_{\sigma(3)}$, which would yield
$a_{\sigma(1)}+a_{\sigma(2)}=a_{\sigma(3)}+a_{\sigma(4)}$, in contradiction with condition (iii).
\item Assume $\alpha=\beta$ and $\car(\K)=2$. Then some $a_i$
is different from $0$ and $\alpha$, which is impossible by case (i) in Corollary \ref{diagcor}.
\item Assume finally that $\alpha \neq \pm \beta$. \\
By cases (iv) and (v) of Corollary \ref{diagcor}, the set
$\calE:=\bigl\{i \in \lcro 1,4\rcro : \; a_i \in \{\alpha,\beta,\alpha+\beta\}\bigr\}$ must have an even cardinal
(because there is an even number of $a_i$'s in $\K^*$ and an even number of $a_i$'s outside of $\{0,\alpha,\beta,\alpha+\beta\}$).
Using the same line of reasoning as in the second point, we see that $\calE$ is not empty (because of condition (iii) and the symmetry condition in cases (iv) and (v) of Corollary \ref{diagcor}).
Hence $\calE$ has two elements, and again, since there are also two of the $a_i$'s outside of $\{0,\alpha,\beta,\alpha+\beta\}$,
their sum is $\alpha+\beta$, so the two elements of $\calE$ cannot be $\alpha$ and $\beta$.
Without loss of generality, we may then assume that $a_1=\alpha$ and $a_2=\alpha+\beta$,
with $a_3$ and $a_4$ outside of $\{0,\alpha,\beta,\alpha+\beta\}$.
Again, cases (iv) and (v) of Corollary \ref{diagcor} would show that $a_3+a_4=\alpha+\beta=a_2$, in contradiction with condition (ii).
\end{itemize}
Finally, there actually exists a quadruple $(a_1,a_2,a_3,a_4) \in \K^4$ which satisfies condition (i) to (iii): indeed, the polynomial
$$P:=\underset{1\leq k<\ell \leq 4}{\prod}(X_k^2-X_\ell^2)\,
\underset{(k,\ell,m)\in \lcro 1,4\rcro^3}{\prod}(X_k+X_\ell-X_m)
\,\underset{\sigma \in \frak{S}_4}{\prod}(X_{\sigma(1)}+X_{\sigma(2)}-X_{\sigma(3)}-X_{\sigma(4)})$$
does not totally vanish on $\K^4$ because $P \neq 0$ and $\K$ is infinite.
Hence there exists a matrix of $\Mat_n(\K)$ which is not a LC of two idempotents, which proves $\ell_n(\K) \geq 3$.
\end{enumerate}

\begin{Rem}
Some of the results on the inability to express matrices with irreducible characteristic polynomials as linear combinations of two idempotents can also be derived from the fact that a simple algebra generated by two non-commuting idempotents over a field $\K$ must be isomorphic to the algebra of 2x2 matrices over a finite extension of $\K$ (see \cite{Laffey}).
\end{Rem}

\section{A review of cyclic matrices, and the key lemma}

\noindent
Let $A \in \Mat_n(\K)$. We say that $A$ is \defterm{cyclic} when
$A \sim C(P)$ for some polynomial $P$ (and then $P=\chi_A$).
A \textbf{good cyclic} matrix is a matrix of the form
$$A=\begin{bmatrix}
a_{1,1} & a_{1,2}  & &  & a_{1,n} \\
1 & a_{2,2} & &   &  \\
0 & \ddots & \ddots & & \vdots \\
\vdots & & & a_{n-1,n-1} & a_{n-1,n} \\
0 & & &  1 & a_{n,n}
\end{bmatrix}$$
with no condition on the $a_{i,j}$'s for $j \geq i$. \\
It is folklore that such a matrix is always cyclic, and, more precisely, that
there exists an upper triangular matrix $T \in \Mat_n(\K)$ with diagonal coefficients all equal to $1$ such that
$T\,A\,T^{-1}=C(\chi_A)$ (this can be seen by performing elementary row and column operations on $A$).

\vskip 2mm
\noindent The following lemma will be the last key to theorem \ref{maindSP}:

\begin{lemme}[Choice of polynomial lemma]\label{cyclicfit}
Let $A \in \Mat_n(\K)$ and $B \in \Mat_p(\K)$ denote two good cyclic matrices, and
$P$ denote a monic polynomial of degree $n+p$ such that $\tr P=\tr A+\tr B$. \\
Then there exists a matrix $D \in \Mat_{n,p}(\K)$ such that
$$\begin{bmatrix}
A & D \\
H_{p,n} & B
\end{bmatrix} \sim C(P).$$
\end{lemme}

\begin{Rem}
The condition on $\tr P$ cannot be done away with since the trace of
$\begin{bmatrix}
A & D \\
H_{p,n} & B
\end{bmatrix}$ is $\tr A+\tr B$.
\end{Rem}

\begin{proof}
We set $M(D):=\begin{bmatrix}
A & D \\
H_{p,n} & B
\end{bmatrix}$. Notice first that $M(D)$ is a good cyclic matrix whatever the choice of $D$,
hence it suffices to show that $D$ can be carefully chosen so that $\chi_{M(D)}=P$. \\
Also, we can replace $A$ and $B$ respectively with $C(\chi_A)$ and $C(\chi_B)$:
indeed, should there be a matrix $D \in \Mat_{n,p}(\K)$ such that
$\begin{bmatrix}
C(\chi_A) & D \\
H_{p,n} & C(\chi_B)
\end{bmatrix}$ has characteristic polynomial $P$,
then there would be two upper triangular matrices $T \in \GL_n(\K)$ and $T' \in \GL_p(\K)$, with diagonal coefficients
all equal to $1$, such that $T \,C(\chi_A)\, T^{-1}=A$ and $T'\, C(\chi_B)\, (T')^{-1}=B$; setting
$T_1:=\begin{bmatrix}
T & 0 \\
0 & T'
\end{bmatrix}$, straightforward computation would then yield
$$T_1\,\begin{bmatrix}
C(\chi_A) & D \\
H_{p,n} & C(\chi_B)
\end{bmatrix}\,T_1^{-1}=M(T\,D\,(T')^{-1}),$$
hence the matrix $TD(T')^{-1}$ would have the required properties. \\
Therefore, we will assume from now on that $A$ and $B$ are respectively the companion matrices of
polynomials $Q=X^n-\underset{k=0}{\overset{n-1}{\sum}}a_k\,X^k$
and $R=X^p-\underset{k=0}{\overset{p-1}{\sum}}b_k\,X^k$. \\
Hence
$$M(D)-X.I_{n+p}=\begin{bmatrix}
-X & 0 & \cdots & 0 & a_0 & d_{1,1} &  & \cdots &  & d_{1,p} \\
1 & -X & 0 & & a_1 &  &  &  &    \\
0 & & \ddots & & \vdots & \vdots & & & & \vdots \\
  & \ddots &  &  -X & a_{n-2} &  \\
  &  & 0 & 1 & -X+a_{n-1} & d_{n,1} & & \cdots & & d_{n,p} \\
  & & & 0 &     1 & -X    & 0 & & 0  & b_0 \\
  &  &  &  & 0 & 1 & -X    & &   & b_1 \\
  &  &  &  &   & 0 & \ddots & \ddots & & \vdots \\
    &  &  &  &   &  & \ddots & 1 & -X & b_{p-2} \\
0  &  &  &  &   &  &  & 0 & 1 & -X+b_{p-1}
\end{bmatrix}.$$
Applying the row operations $L_i \leftarrow L_i+XL_{i+1}$ for
$i$ downward from $n-1$ to $1$, we obtain that $M(D)-X.I_{n+p}$ has the same determinant has
$$\begin{bmatrix}
0 & 0 & \cdots & 0 & -Q(X) & P_1(X) &  & \cdots &  & P_p(X) \\
1 & 0 & 0 & & a_1 & ? &  &  & & ?  \\
0 & & \ddots & & \vdots & \vdots & & & & \vdots \\
  & \ddots &  &  0 & ? &  \\
  &  & 0 & 1 & -X+a_{n-1} & ? & & \cdots & & ? \\
  & & & 0 &     1 & -X    & 0 & & 0  & b_0 \\
  &  &  &  & 0 & 1 & -X    & &   & b_1 \\
  &  &  &  &   & 0 & \ddots & \ddots & & \vdots \\
    &  &  &  &   &  & \ddots & 1 & -X & b_{p-2} \\
0  &  &  &  &   &  &  & 0 & 1 & -X+b_{p-1}
\end{bmatrix},$$
where, for all $j \in \lcro 1,p\rcro$, $P_j:=\underset{k=0}{\overset{n-1}{\sum}}d_{k+1,j}\,X^j$. \\
By developing inductively this determinant along the first column, we get:
$$\det(M(D)-X.I_{n+p})
=(-1)^{n-1}
\det\begin{bmatrix}
-Q(X) & P_1(X) &  & \cdots &  & P_p(X) \\
1 & -X    & 0 & & 0  & b_0 \\
0 & 1 & -X    & &   & b_1 \\
  & 0 & \ddots & \ddots & & \vdots \\
   &  & \ddots & 1 & -X & b_{p-2} \\
0    &  &  & 0 & 1 & -X+b_{p-1}
\end{bmatrix}.$$
Development of this last determinant along the first row finally yields:
$$\chi_{M(D)}=Q(X)\,R(X)-\underset{j=1}{\overset{p}{\sum}}\,P_j(X)\, R_{p-j}(X)$$
where, for $j \in \lcro 0,p-1\rcro$, we have set $R_j(X):=X^j-\underset{k=0}{\overset{j-1}{\sum}}b_{k+p-j}\,X^k$.
Proving that there is a $D \in \Mat_{n,p}(\K)$ such that
$\chi_{M(D)}=P$ is thus equivalent to proving that there are
$p$ polynomials $P_1,\dots,P_p$ in $\K_{n-1}[X]$ (i.e. of degree at most $n-1$) such that
$$P-Q\,R=\underset{j=1}{\overset{p}{\sum}}P_j\, R_{p-j}.$$
This however comes readily by noticing that the condition on the degree of $P$ and its trace show that
$\deg(P-Q\,R)<n+p-1$ and that the $(n+p-1)$-tuple
$(R_0,R_1,\dots,R_{p-2},R_{p-1},X\,R_{p-1},X^2\,R_{p-1},\dots,X^{n-1}\,R_{p-1})$ is a basis of
$\K_{n+p-2}[X]$ (since it features $n+p-1$ polynomials, with one of degree $k$ for every $k \in \lcro 0,n+p-2\rcro$).
\end{proof}

\noindent Finally, this basic lemma of reduction theory will be used
at crucial steps in this paper:

\begin{lemme}
Let $A \in \Mat_n(\K)$, $B \in \Mat_p(\K)$, and $C \in \Mat_{n,p}(\K)$.
Assume $\chi_A$ and $\chi_B$ are mutually prime. Then
$$\begin{bmatrix}
A & C \\
0 & B
\end{bmatrix} \sim \begin{bmatrix}
A & 0 \\
0 & B
\end{bmatrix}$$
\end{lemme}

\begin{Rem}
This is a special case of Roth's theorem \cite{Roth}. For alternative proofs and extensions, see
\cite{Guralnick1} and \cite{Guralnick2}. 
\end{Rem}

\begin{proof}
For any $M \in \Mat_{n,p}(\K)$, we have:
$$\begin{bmatrix}
I_n & M \\
0 & I_p
\end{bmatrix}
\begin{bmatrix}
A & C \\
0 & B
\end{bmatrix}
\begin{bmatrix}
I_n & M \\
0 & I_p
\end{bmatrix}^{-1}=\begin{bmatrix}
A & C+MB-AM \\
0 & B
\end{bmatrix}.$$
It thus suffices to prove that the endomorphism $$\begin{cases}
\Mat_{n,p}(\K) & \longrightarrow \Mat_{n,p}(\K) \\
M & \longmapsto AM-MB
\end{cases}$$
is onto, which is true if it is one-to-one. Let $M \in \Mat_{n,p}(\K)$
such that $AM=MB$. Then the matrix
$\begin{bmatrix}
I_n & M \\
0 & I_p
\end{bmatrix}$ commutes with
$\begin{bmatrix}
A & 0 \\
0 & B
\end{bmatrix}$. Since $A$ and $B$ have mutually prime annihilator polynomials,
this forces $\begin{bmatrix}
I_n & M \\
0 & I_p
\end{bmatrix}$ to stabilize $\{0\} \times \K^p$
(seen as a linear subspace of $\K^{n+p}$), hence $M=0$, which completes the proof.
\end{proof}

\section{Every matrix is a linear combinations of three idempotents}\label{LC3}

In this section, we fix a matrix $A \in \Mat_n(\K)$ and prove that it can be decomposed
as an LC of three idempotents. This will complete the proof of Theorem \ref{maindSP}.
The basic idea is to add $A$ to a scalar multiple of
an idempotent in order to obtain a linear combination of two idempotents. \\
In the course of the proof, we will use the following basic fact repeatedly
(cf. \cite{Gantmacher}):
when $P$ and $Q$ denote two monic polynomials which are mutually prime, one has
$$C(P\,Q) \sim \begin{bmatrix}
C(P) & 0 \\
0 & C(Q)
\end{bmatrix}.$$
Using this and a rational canonical form, we see that any matrix is similar to
$D\bigl(C(P_1),\dots,C(P_N)\bigr)$, where $P_1,\dots,P_N$ are monic polynomials each of which
has essentially one irreducible divisor (this is the primary canonical form for the matrix).

\subsection{When the minimal polynomial of $A$ is a power of an irreducible polynomial}

Here, we assume that the minimal polynomial of $A$ is a power of an irreducible monic polynomial
$P=X^p-\underset{k=0}{\overset{p-1}{\sum}}a_k X^k$.
If $p=1$, then there is some $\alpha \in \K$ and some nilpotent matrix $N$ such that
$A=\alpha.I_n+N$, so $A$ is $(\alpha,1,-1)$-composite.
Assume now that $p \geq 2$. \\
For any $k \in \N^*$, set
$$M_k:=\begin{bmatrix}
C(P) & 0 & \dots & 0 \\
H_{p,p} & C(P) &  & \vdots \\
 & \ddots & \ddots & 0\\
0 & &  H_{p,p} & C(P)
\end{bmatrix} \in \Mat_{kp}(\K).$$
By the generalized Jordan reduction theorem, there are integers
$k_1,\dots,k_N$ such that
$$A \sim D(M_{k_1},\dots,M_{k_N}),$$
so we lose no generality assuming $A=D(M_{k_1},\dots,M_{k_N})$.
\begin{itemize}
\item \emph{The case $\tr(P) \neq 0$.} \\
Set then $\alpha:=\dfrac{1}{\tr(P)}$,
$$G:=\begin{bmatrix}
0 &   & 0 & \alpha\, a_0 \\
\vdots &  \ddots & \vdots & \vdots \\
0 &  & 0  & \alpha\, a_{p-2} \\
0 & \cdots & 0 & 1 \\
\end{bmatrix} \in \Mat_p(\K)$$ and, for all $k \in \N^*$,
$G_k:=D(G,\dots,G) \in \Mat_{kp}(\K)$.
The matrix $B:=D(G_{k_1},\dots,G_{k_N})$ is clearly idempotent, whilst
$A-\frac{1}{\alpha}\,B$ is clearly nilpotent, hence $A$ is a $\bigl(\tr(P),1,-1\bigr)$-composite.

\vskip 2mm
\item \emph{The case $\tr P=0$.} \\
Set now
$$G:=\begin{bmatrix}
0 &   & 0 & a_0 \\
\vdots &  \ddots & \vdots & \vdots \\
0 &  & 0  & a_{p-2} \\
0 & \cdots & 0 & 1 \\
\end{bmatrix} \in \Mat_p(\K)$$ and, for all $k \in \N^*$,
$$G_k:=\begin{bmatrix}
G & 0 & \dots & 0 \\
H_{p,p} & G & &  \vdots \\
 & \ddots & \ddots & 0 \\
0 & & H_{p,p} & G
\end{bmatrix} \in \Mat_{kp}(\K).$$
Again, the matrix $B:=D(G_{k_1},\dots,G_{k_N})$
is idempotent, and this time
$$A-B \sim D\bigl(C(X^p+X^{p-1}),\dots,C(X^p+X^{p-1})\bigr) \sim D\bigl(C(X^{p-1}),\dots,C(X^{p-1}),-I_q\bigr)$$
for some integer $q$. It follows that $A-B$ is a difference of two idempotents,
hence $A$ is a $(1,1,-1)$-composite.
\end{itemize}

\subsection{When the minimal polynomial of $A$ is not a power of an irreducible polynomial}

We now assume that the minimal polynomial of $A$ has two different monic
irreducible divisors. \\
We will first prove the following fact:

\begin{lemme}
Assume the  minimal polynomial of $A$ has two different monic
irreducible divisors. Then there are two distinct $\alpha$ and $\beta$
in $\K$, integers $p$ and $q$ (possibly zero), non constant monic polynomials $P_1,\dots,P_r$
and $Q_1,\dots,Q_s$ (with $r \geq 1$ and $s \geq 1$)
such that $\deg P_j \geq 2$ for all $j \geq 2$, $\deg Q_k \geq 2$ for all $k \leq s-1$,
at most one of the polynomials $P_1$ and $Q_s$ has degree $1$, any $P_i$ is prime to any $Q_j$, and
$$A \sim D\bigl(\alpha.I_p,\beta.I_q,C(P_1),\dots,C(P_r),C(Q_1),\dots,C(Q_s)\bigr)$$
\end{lemme}

\begin{proof}
We start by reducing $A$ to a primary canonical form, so $A$ is similar to a block-diagonal matrix
of the form
$$A'=D\bigl(\alpha_1.I_{n_1},\alpha_2.I_{n_2},\dots,\alpha_N.I_{n_N},C(P_1^{a_1}),\dots,C(P_m^{a_m})\bigr)$$
where $P_1,\dots,P_m$ are irreducible monic polynomials of degree greater or equal to $2$, and
$n_1 \geq n_2 \geq \dots \geq n_N$ (possibly with $n_1=0$ or $n_2=0$, for sake of generality).
We immediately leave aside the trivial case where $N \leq 3$ and $m=0$.

\begin{itemize}
\item If $n_1>0$ and $n_2=0$, then we immediately obtain a similarity
$$A \sim D\bigl(\alpha_1.I_{n_1-1},\alpha_1,C(Q_1),\dots,C(Q_q),C(R_1),\dots,C(R_r)\bigr),$$
where $r \geq 1$, $q \geq 0$, the $Q_k$'s are powers of $X-\alpha_1$ with $\deg Q_k \geq 2$, the
$R_k$'s have degree greater or equal to $2$ and $\alpha_1$ is not a root of any of them.

\item If $n_2>0$, $n_3=0$, and $\alpha_1$ is a root of some $P_i$, then
we obtain a similarity
$$A \sim D\bigl(\alpha_1.I_{n_1},\alpha_2.I_{n_2-1},C(Q_1),\dots,C(Q_q),C(R_1),\dots,C(R_r),\alpha_2\bigr),$$
with $q \geq 1$, $r \geq 0$, whilst the $Q_k$'s and the $R_k$'s have the same properties as
in the first point.

\item If $n_2>0$, $n_3=0$, and $\alpha_1$ is a root of none of the $P_i$'s, then
we have a similarity
$$A \sim D\bigl(\alpha_1.I_{n_1-1},\alpha_2.I_{n_2},\alpha_1,C(P_1),\dots,C(P_m)\bigr),$$
and $m \geq 1$.

\item Finally, if $n_3>0$, then we can use the similarities
$D(\alpha_2,\dots,\alpha_j) \sim C\bigl((X-\alpha_2)\cdots (X-\alpha_j)\bigr)$ for $j \in \lcro 3,N\rcro$
to obtain a similarity
$$A \sim D\bigl(\alpha_1.I_{n_1-1},\alpha_2.I_{n_2-n_3},\alpha_1,C(Q_1),\dots,C(Q_q),C(R_1),\dots,C(R_r)\bigr),$$
where $q \geq 0$, $r \geq 1$, and the $Q_k$'s and $R_j$'s have the same properties as in the first point.
\end{itemize}
In any case, the lemma is proven.
\end{proof}

\noindent We now set $\alpha$, $\beta$, $p$, $q$, $P_1,\dots,P_r$ and $Q_1,\dots,Q_s$ as in the above lemma, so
$$A \sim D\bigl(\alpha.I_p,\beta.I_q,C(P_1),\dots,C(P_r),C(Q_1),\dots,C(Q_s)\bigr).$$
We will now focus on the block-diagonal matrix
$$B:=D\bigl(C(P_1),\dots,C(P_r),C(Q_1),\dots,C(Q_s)\bigr).$$
We let $t$ denote the size of $B$. Our next aim is the following key lemma:

\begin{lemme}
Let $P \in \K[X]$ be a monic polynomial of degree $t$ such that $\tr P \neq \tr B$.
Then there exists an idempotent $Q \in \Mat_t(\K)$ and a scalar $\delta$ such that
$$B-\delta\,Q \,\sim\, C(P).$$
\end{lemme}

\begin{proof}
For $i \in \lcro 1,r\rcro$ and $j \in \lcro 1,s\rcro$, set
$n_i:=\deg P_i$ and $m_j:=\deg Q_j$. \\
Define $\lambda:=(r+s-1).1_\K$ if $(r+s-1).1_\K \neq 0$, or else
$\lambda:=(r+s).1_\K$ (so that $\lambda \neq 0$ in any case). \\
For $k \in \N^*$, recall that $F_k=D(0,\dots,0,1) \in \Mat_k(\K)$. \\
Let $\delta \in \K^*$ and define $R(\delta)$ as:
$$\begin{bmatrix}
F_{n_1} & 0 &  & & & & & & 0\\
-\frac{1}{\delta}\,H_{n_2,n_1} & F_{n_2} & \ddots & & \\
0 & \ddots & \ddots & 0  \\
 &  & -\frac{1}{\delta}\,H_{n_r,n_{r-1}} & F_{n_r}& 0 \\
 & & & -\frac{1}{\delta}\, H_{m_1,n_r} & F_{m_1} & 0 & & & \\
\vdots & & & & -\frac{1}{\delta}\,H_{m_2,m_1} & F_{m_2} & & & \\
 & & & &   & \ddots & \ddots & \ddots & \\
 & & & &   &   & -\frac{1}{\delta}\,H_{m_{s-1},m_{s-2}} & F_{m_{s-1}} & \\
0 & & \cdots & & \cdots & & 0 & -\frac{1}{\delta}\,H_{m_s,m_{s-1}} & 0
\end{bmatrix}.$$
If $(r+s-1).1_\K \neq 0$, set
$$Q(\delta):=R(\delta).$$
If $(r+s-1).1_\K =0$ and $n_1>1$, set
$$Q(\delta):=D(1,0,\dots,0)+R(\delta).$$
If $(r+s-1).1_\K =0$ and $n_1=1$, then $m_s > 1$ and we can therefore set
$$Q(\delta):=D(0,\dots,0,1)+R(\delta).$$
In any case:
\begin{itemize}
\item $Q(\delta)$ is idempotent;
\item $\tr Q(\delta)=\lambda$;
\item There are good cyclic matrices $B'_1$ and $B'_2$ such that:
$$B-\delta.Q(\delta)=\begin{bmatrix}
B'_1 & 0 \\
H_{M,N} & B'_2
\end{bmatrix},$$
where $M=\underset{k=1}{\overset{s}{\sum}}m_k$ and $N=\underset{k=1}{\overset{r}{\sum}}n_k$.
\end{itemize}
We now choose $\delta:=\dfrac{\tr B-\tr P}{\lambda}$, so that
$$\tr B'_1+\tr B'_2=\tr(B-\delta.Q(\delta))=\tr B-\lambda\,\delta=\tr P.$$
By Lemma \ref{cyclicfit}, there exists a matrix $D \in \Mat_{N,M}(\K)$ such that
$$\begin{bmatrix}
B'_1 & D \\
H_{M,N} & B'_2
\end{bmatrix} \,\sim\, C(P).$$
Let us now decompose
$$B=\begin{bmatrix}
B_1 & 0 \\
0 & B_2
\end{bmatrix} \quad \text{with $B_1 \in \Mat_N(\K)$ and $B_2 \in \Mat_M(\K)$.}$$
Notice that the assumptions on the polynomials $P_i$ and $Q_k$ imply that $\chi_{B_1}$ and $\chi_{B_2}$ are mutually prime, so
$$B \sim B':=\begin{bmatrix}
B_1 & D \\
0 & B_2
\end{bmatrix}.$$
However,
$$B'-\delta.Q(\delta)=\begin{bmatrix}
B'_1 & D \\
H_{M,N} & B'_2
\end{bmatrix} \,\sim\, C(P)$$
so there exists an idempotent $Q'$ similar to $Q(\delta)$
with $B-\delta.Q' \sim C(P)$.
\end{proof}

\vskip 2mm
\noindent We can now complete our proof.
Let $P \in \K[X]$ denote a monic polynomial of degree $t$ such that $\tr P \neq \tr B$.
Then there exists an idempotent $Q'$ and a scalar $\delta$ such that
$B-\delta.Q' \sim C(P)$, so
$Q'':=\begin{bmatrix}
0 & 0 \\
0 & Q
\end{bmatrix} \in \Mat_n(\K)$ is also an idempotent and
$A-\delta.Q'' \sim \begin{bmatrix}
\alpha.I_r & 0 & 0 \\
0 & \beta.I_s & 0 \\
0 & 0 & C(P)
\end{bmatrix}$.
The proof of Theorem \ref{maindSP} will then be completed should we establish the following lemma:

\begin{lemme}
Let $(\alpha,\beta)\in \K^2$ such that $\alpha \neq \beta$, and
$(r,s,t)\in \N^2 \times \N^*$. Let $\gamma \in \K$. \\
Then there exists a monic polynomial
$P \in \K[X]$ of degree $t$ such that $\tr P \neq \gamma$ and
the block-diagonal matrix
$M(P):=\begin{bmatrix}
\alpha.I_r & 0 & 0 \\
0 & \beta.I_s & 0 \\
0 & 0 & C(P)
\end{bmatrix}$
is a linear combination of two idempotents.
\end{lemme}

\begin{proof}
\vskip 2mm
\item \emph{The case $\alpha =0$ or $\beta=0$.} \\
Without loss of generality, we may actually assume $\beta=0$. \\
The two polynomials $P_1=X^t$ and $P_2=X^{t-1}(X-\alpha)$
have then different traces (one of which is different from $\gamma$)
with
$$M(P_1) \sim D\bigl(\alpha.I_r,0.I_s,C(X^t)\bigr) \quad \text{and} \quad
M(P_2) \sim D\bigl(\alpha.I_{r+1},0.I_s,C(X^{t-1})\bigr).$$
Corollary \ref{nilpotent} then shows that $M(P_1)$ and $M(P_2)$ are $(\alpha,-\alpha)$-composites,
hence one of the polynomials $P_1$ or $P_2$ is a solution to our problem.

\vskip 2mm
\item \emph{The case $\alpha \neq 0$ and $\beta \neq 0$.}
\begin{itemize}
\item Assume $t=2\,t'$ for some $t' \in \N$. Then
the polynomials $P_1=(X-\alpha)^{t'}(X-\beta)^{t'}$ and
$P_2=(X-\alpha)^{t'}(X-\beta)^{t'-1}(X-\alpha-\beta)$
have distinct traces. Also
$$M(P_1) \sim D\bigl(\alpha.I_r,\beta.I_s,C((X-\alpha)^{t'}),C((X-\beta)^{t'})\bigr)$$
and
$$M(P_2) \sim D\bigl(\alpha+\beta,\alpha.I_r,\beta.I_s,C((X-\alpha)^{t'}),C((X-\beta)^{t'-1})\bigr),$$
so Corollary \ref{alphaneqbeta} shows that both matrices $M(P_1)$ and $M(P_2)$ are
$(\alpha,\beta)$-composites.
\item Assume $t=2\,t'+1$ for some integer $t'$.
Then
the polynomials $P_1=(X-\alpha)^{t'+1}(X-\beta)^{t'}$ and
$P_2=(X-\alpha)^{t'}(X-\beta)^{t'}(X-\alpha-\beta)$
have distinct traces and
$$M(P_1) \sim D\bigl(\alpha.I_r,\beta.I_s,C((X-\alpha)^{t'+1}),C((X-\beta)^{t'}\bigr)$$
and
$$M(P_2) \sim D\bigl(\alpha+\beta,\alpha.I_r,\beta.I_s,C((X-\alpha)^{t'}),C((X-\beta)^{t'})\bigr)$$
(with the convention that $C(1)$ is the zero matrix of $\Mat_0(\K)$), and
again both $M(P_1)$ and $M(P_2)$ are $(\alpha,\beta)$-composites.
\end{itemize}
\end{proof}

\end{document}